\documentclass[12pt]{article}
\usepackage[latin1]{inputenc}

\usepackage{amssymb}
\usepackage{amstext}
\usepackage{amsmath}
\usepackage{amsfonts}
\usepackage{amssymb}

\usepackage{authblk}

\usepackage{graphicx}
\usepackage{multicol}        
\usepackage{verbatim}
\usepackage{epstopdf}

\newcommand{\NORM}[1]{\left \| {#1} \right \| }
\newcommand{\QT}[1] {"{#1}"}
\newcommand{\MIN}[1] {\mathop{\mbox{min}}_{#1}}
\newcommand{\MAX}[1] {\mathop{\mbox{max}}_{#1}}

\newcommand{\ARGMIN}[1] {\mathop{\mbox{argmin}}_{#1}}

\newcommand{\E}[2]{ \mathcal{E} \! \left. \left [ \, {#1} \, \right |  {#2} \, \right ] }                        

\newcommand{\V}[1]{#1} 

\begin{document}

\title{Approximate Dynamic Programming based on High Dimensional Model Representation \thanks{This work has been funded by the GA\v{C}R project P 102/11/0437}}

\author{Miroslav Pi\v{s}t\v{e}k\thanks{email: miroslav.pistek@gmail.com} \\Institute of Information Theory and Automation, Academy of
Sciences of the Czech Republic, Pod~vod\'{a}renskou v\v{e}\v{z}\'{i} 4, CZ~18208 Praha~8, Czech Republic}

\maketitle 

\begin{abstract}
This article introduces an algorithm for implicit High Dimensional Model Representation (HDMR) of the Bellman equation. This approximation technique reduces memory demands of the algorithm considerably. Moreover, we show that HDMR enables fast approximate minimization which is essential for evaluation of the Bellman function. In each time step, the problem of parametrized HDMR minimization is relaxed into trust region problems, all sharing the same matrix. Finding its eigenvalue decomposition, we effectively achieve estimates of all minima. Their full-domain representation is avoided by HDMR and then the same approach is used recursively in the next time step. An illustrative example of N-armed bandid problem is included. We assume that the newly established connection between approximate HDMR minimization and the trust region problem can be beneficial also to many other applications. 


\end{abstract}

\maketitle 

\section{Introduction}
\label{Section.Introduction}

The main focus of this article is to develop an approximate tool
suitable for enlarging the class of computationally feasible
decision-making problems. It copes with the principal problem within the stochastic dynamic programming,
which is known as the {\it curse of dimensionality}, see  \cite{Powell2007_ApproximateDP}.  The central notion of stochastic dynamic programming is the Bellman function, see for instance \cite{Kusner1971_Stochastic_Control}. Once we are able to find and store this function, it is easy to derive the optimal strategy. However, the exact calculation of the Bellman function is computationally infeasible in the majority of practical applications, and also its representation as a lookup-table is intractable. 

Next, we present a survey of approximate solutions to these problems. One way to reduce the size of the lookup-table is to aggregate the state space of the original problem into smaller sets. As it is not clear how to pick the best level of aggregation, several methods of multiple-level aggregation are developed \cite{Powell2008_MultipleAggregation}. A similar way to lookup-table reduction is approximation of the Bellman function which does not require any simplifications in the state space. A grid-based approximation with value interpolation is a typical example of such method \cite{Hauskrecht2000_Approximations_for_POMDP}. The Bellman function can also be estimated using regression models which are able to exploit specialized structures (\QT{basis functions}) in the state space \cite{LeBlanc1996_Regression}. Nonetheless, such methods are suitable for maximally hundreds of regression parameters.  Another tool suitable for approximation is the artificial neural networks utilized to learn the shape of the Bellman function, see \cite{Miller1995_Neural} and references therein. Based on random sampling of the state space, a variety of Monte Carlo methods may be also applied, see for instance \cite{Luus2000_Iterative}. Temporal Difference methods are of quite a different nature. Opposite to the algorithm developed later, they do not operate with system model. They use simulated or experience-based sampling of system trajectories instead, and thus they have no ambition to cover the whole state space. Nonetheless, they definitely do well for many real-world problems  \cite{Gosavi2009_RL_Survey,Jaakkola1995_RL_for_POMDP,Sutton1998_Reinforcement_Learning}. 
 
In this article, we develop a new approximate technique which considerably reduces both computational and memory demands of a decision-making problem. To this end, an approximation tool called High Dimensional Model Representations (thereinafter \QT{HDMR}) is useful  \cite{Rabitz1999_HDMR_foundations}. It was applied to continuous function approximation in calculating reliability of uncertain mechanical systems \cite{Rahman2008_polynomial_decomposition}. It was also utilized for solution of stochastic partial differential equations \cite{Zabaras2010_adaptive_HDMR} and compared to Monte Carlo sampling. Another application of HDMR was volatility calibration \cite{Kucherenko2009_volatility_calibration} where it was compared to cubic spline approximation.  These successful implementations of HDMR in other fields encourage us to apply it to approximate dynamic programming.  In the previous applications, it was used mainly for reducing the amount of data. The memory space necessary to store all the values of function $g(x_1,\ldots,x_d)$ grows
exponentially with the dimension $d$, whereas the size growth of HDMR components is just quadratic in $d$. This is, of course, important even in our case, 
but the newly established fact that HDMR permits fast approximate minimization may be even more essential in the context of the decision making theory. 

The outline of this work is as follows. Section \ref{Section.HDMR} deals with the approximation
technique of HDMR, which is determined by a system of linear equations. Its linearity does not match with the inherently non-linear Bellman equation. On that account, an algorithm for approximate minimization of function having HDMR form is developed in Section \ref{Section.HDMRminimization}. Then, the current state of the art in the decision making theory is summarized at the beginning of Section \ref{Sec.ADPHDMR}.  Next, a viable technique for approximate decision making based on HDMR is introduced there, and then the $N$-armed bandit problem is tackled as an example. Section
\ref{Section.Conclusion} is devoted to conclusion.

Throughout this work, a few general conventions are followed. The domain of the quantity $x$ is denoted $X$, $x\in X$, $|X|$ denotes the count of elements of finite set $X$. Next, $x_m$ denotes $m$-th coordinate for vector valued quantity $\V{x} \subset \mathbb{R}^d$, $\V{x} =
(x_1,\ldots,x_d)$. This convention holds with one exception: if we use letter $t$ as a subscript, e.g. $\V{x}_t$, it stands for quantity $\V{x}$ at the time instant $t \in T$ with $T$ finite. Next, we
reserve letter \QT{$f$} for conditional probability density functions, arguments in the condition are separated by \QT{$|$} in the argument list. For the domain of function $h(x)$ we use $\mbox{dom}(h)$, and HDMR of $h(x)$ is marked by $\tilde{h}(x)$ with several  exceptions pointed out later.
\section{High Dimensional Model Representation}
\label{Section.HDMR}
The approximation technique of HDMR has a particularly simple form. For a general function $g(\V{x})$, the second order HDMR $\tilde{g}(\V{x})$ reads
\begin{eqnarray}
\label{HDMR:Example}
g(\V{x}) & \approx & \tilde{g}(\V{x})  =
\tilde{g}(x_{1},x_{2},\ldots,x_{d})  =  \\
&& \tilde{g}_{\emptyset}+ \sum_{m=1}^{d} \tilde{g}_{m}(x_{m})+
\sum_{m=1}^{d-1} \sum_{n=m+1}^{d}\tilde{g}_{mn}(x_{m},x_{n}) \nonumber .
\end{eqnarray} 
Here, 
$\tilde{g}_{\emptyset}$ denotes a constant value over $\mbox{dom}(g)$; one-dimensional functions $\tilde{g}_{m}(x_{m})$ describe independent effects of each particular coordinate $x_{m}$, and two-dimensional functions $\tilde{g}_{mn}(x_{m},x_{n})$ represent the
joint effect of coordinates $x_{m}$ and $x_{n}$. In the context of HDMR, these functions are called zero-order, first-order, and second-order components of HDMR, respectively.
Experience shows that second-order HDMR provides a sufficient approximation
of $g(\V{x})$ as only low-order correlations amongst the input variables have a significant impact upon the outputs of a typical model \cite{Rahman2008_polynomial_decomposition,Zabaras2010_adaptive_HDMR,Kucherenko2009_volatility_calibration}. 

There are many ways how to construct HDMR \cite{Rabitz1999_HDMR_foundations,Dermiralp2003_HDMR_varieties}.
To reduce this ambiguity, it is thus necessary to
formalize its desired properties. The function Hilbert space $L^2(X)$ is a useful concept for the
function approximation. It is a space of real functions
defined over a set $X$ with the finite norm $\NORM g$ defined as follows
\begin{equation}
\label{NormDef}
\NORM{g}^2 := \int_{X} g(\V{x})^2 \,d\V{x}.
\end{equation} 
Then, the optimal HDMR of the function $g \in L^2(X)$
is defined as a minimizer of the approximation error $\NORM{g- \tilde{g}}$. The uniqueness of projection on closed subspaces of $L^2(X)$ implies the uniqueness of minimizing function $\tilde{g}(\V{x})$ matching this form
\begin{equation} \label{GTILDE}
\tilde{g}(\V{x})  = \tilde{g}_{\emptyset}+ \sum_{m=1}^d \tilde{g}_{m}(x_{m})+ \frac{1}{2} 
\sum_{m,n=1}^{d} \tilde{g}_{mn}(x_{m},x_{n}),
\end{equation}
where we slightly generalized (\ref{HDMR:Example}) to better fit our needs. Nonetheless, there may exist various components
$\tilde{g}_{\emptyset}$, $\tilde{g}_{m}$ and $\tilde{g}_{mn}$ summing up to the same $\tilde{g}$. 

Now, let $X$ be $d$-dimensional product of finite sets $X_i$
\begin{equation} \label{d-dim_X}
X = \prod_{i=1}^d X_i,
\end{equation} 
and let the integration in (\ref{NormDef}) be summation over $X$. Next, for any subset of indices  $I \subset \{ 1,\ldots, d \}$ we define 
\begin{equation}
X^\bot_I := \prod_{ i \in \{ 1, \ldots, d \} \setminus I} X_i.
\end{equation}
Then, the optimal HDMR of $\tilde{g}$ may be obtained from marginal operators defined for function $g$ as
\begin{eqnarray} \label{HDMR.marginals}
\mbox{M}_{\emptyset}[g] & := & \sum_{y \in X} g(y_1,\ldots,y_{d})  \\
\mbox{M}_{m}[g](x_m) & := & \sum_{y \in X_m^\bot} g(y_1,\ldots,y_{m-1},x_m,y_{m+1},\ldots,y_{d})  \nonumber \\
\mbox{M}_{mn}[g](x_m,x_n) & := & \sum_{y \in X_{mn}^\bot} \! g(y_1,\ldots,y_{m-1},x_m,y_{m+1},\ldots, x_{n},\ldots,y_{d}) \nonumber.
\end{eqnarray}
The formulae for HDMR components of the optimal $\tilde{g}(\V{x})$ read 
\begin{eqnarray} \label{HDMR.approximation}
\tilde{g}_{\emptyset} & := & \frac{1}{|X|} \, \mbox{M}_{\emptyset}[g] \\
\tilde{g}_{m}(x_{m}) & := & \frac{1}{|X_m^\bot|} \, \mbox{M}_{m}[g](x_m) - \tilde{g}_{\emptyset}  \nonumber \\
\tilde{g}_{mn}(x_m,x_n) & := & \frac{1}{|X_{mn}^\bot|} \, \mbox{M}_{mn}[g](x_m,x_n)  - \tilde{g}_m(x_m)- \tilde{g}_n(x_n) - \tilde{g}_{\emptyset} \nonumber \\ \tilde{g}_{mm}(x_m,x_m) &=& 0 \nonumber .
\end{eqnarray}
The proposed variant of approximation matches \QT{ANOVA-HDMR} in \cite{Rabitz1999_HDMR_foundations}. From equations (\ref{HDMR.approximation}) we observe that  identities 
\begin{eqnarray} \label{HDMR.approximation_mean}
\sum_{x_{m} \in X_m} \tilde{g}_{m}(x_{m})  = 0 \\ \sum_{x_{m} \in X_m} \sum_{x_{n} \in X_n} \tilde{g}_{mn}(x_m,x_n) = 0 \nonumber
\end{eqnarray}
hold for all $m,n \in \{1,\ldots,d\}$. In fact, construction (\ref{HDMR.approximation}) was intentionally designed to satisfy (\ref{HDMR.approximation_mean}) in order to provide uniqueness of all HDMR components \cite{Rabitz1999_HDMR_foundations}. In our setting, however,  identities (\ref{HDMR.approximation_mean}) play also another important role in Section \ref{Section.HDMRminimization}.

Finally, we note that this simple construction of HDMR is beneficial to our application, as the domain of the Bellman function could be too large to operate with all the function values at once. Still, its HDMR components can be computed by pointwise evaluation of the function values which are immediately added to proper sums in (\ref{HDMR.approximation}). Next, we show that such convenient form of HDMR may be constructed even in a more general setting.
\subsection{Weighted HDMR}
\label{weighted_H}
A more difficult construction of HDMR may occur in practise if the
approximated function $g(\V{x})$ is defined only on a strict subset of $X$, $\mbox{dom}(g) = R \subsetneq X$. Or, if the full domain $X$ is too large to handle, and thus we search only for some approximation to the optimal HDMR, which may be constructed from samples of $g(x)$ taken with respect to a smaller set, and so we have $x \in R \subsetneq X$ again. Both these situations may clearly arise in the decision-making theory.

Under such conditions, it is important not to consider points $X \setminus R$ in the computation 
of the approximation error. Thus, instead of (\ref{NormDef}) we have to use a weighted norm
\begin{equation}
\label{WeightedNormDef}
\NORM{g}^2_{\chi_R} := \int_{X} \chi_R(\V{x})\,g(\V{x})^2\,\mbox{d}\V{x} = \int_{R} g(\V{x})^2\,\mbox{d}\V{x}
\end{equation} 
with a weight equal to characteristic function
\begin{equation}
\chi_R(\V{x}) := 1 \quad \mbox{for} \quad \V{x} \in R, \qquad 
\chi_R(\V{x}) := 0 \quad \mbox{for} \quad \V{x} \not \in R.
\label{CharDef}
\end{equation} 
We note that for the case of product weight satisfying $w(\V{x}) = \prod_{i=1}^n w_i(x_i)$, the optimal HDMR with respect to $\NORM{g}^2_{w}$ may be obtained identically to (\ref{HDMR.approximation}), see again \cite{Rabitz1999_HDMR_foundations}.  This is, however, not possible for an intrinsically non-product weight $\chi_R(\V{x})$. 

Yet, we can directly minimize the approximation error with respect to (\ref{WeightedNormDef}), but instead of component-wise decoupled equations  (\ref{HDMR.approximation}), we obtain one large linear system determining all the optimal HDMR components of $\tilde{g}(\V{x})$, see \cite{Pistek2009_BE_Approximation}. For smaller problems this system may be computationally feasible; however, a more convenient way is to slightly redefine our task. Instead of searching for an optimal approximation within the class of all functions having HDMR form (\ref{GTILDE}), we search for it within a smaller class of such HDMR functions that are determined by decoupled formulae as in (\ref{HDMR.approximation}). The crucial property is the mutual independence of HDMR components: $\tilde{g}_{\emptyset}$ does not depend on any other HDMR component, each $\tilde{g}_{m}(x_{m})$ depends only on $\tilde{g}_{\emptyset}$, and finally each $\tilde{g}_{mn}(x_m,x_n)$ depends only on $\tilde{g}_m(x_m)$, $\tilde{g}_n(x_n)$ and $\tilde{g}_{\emptyset}$. By enforcing only these hierarchical relations we obtain an easier computation of HDMR components of $\tilde{g}(\V{x})$ at the price of worse approximation.

We build such second order HDMR in three steps. First, we compute zero order component $\tilde{g}_\emptyset$ in such a way that it minimizes the approximation error $\NORM{g- \tilde{g}_\emptyset}_{\chi_R}$. In the next step we fix this component and find such first order components $\tilde{g}_m(x_m)$ that minimize approximation error $\NORM{g- \tilde{g}_\emptyset - \tilde{g}_m}_{\chi_R}$ with respect to $\tilde{g}_\emptyset$. Finally, we find second order components $\tilde{g}_{mn}(x_m,x_n)$ as minimizers of $\NORM{g- \tilde{g}_\emptyset - \tilde{g}_m  - \tilde{g}_n - \tilde{g}_{mn}}_{\chi_R}$ with $\tilde{g}_\emptyset$, $\tilde{g}_m(x_m)$, and $\tilde{g}_n(x_n)$ kept fixed.
The optimality conditions for such HDMR may be derived in three steps where each step is analogous to the original derivation of the full HDMR \cite{Rabitz1999_HDMR_foundations}. Thus, we obtain the following decoupled sytem of equations determining HDMR components of $\tilde{g}(\V{x})$ 
\begin{eqnarray} \label{wHDMR.approximation}
\tilde{g}_{\emptyset} & := & \frac{\mbox{M}_{\emptyset}[ \chi_R. g]}{\mbox{M}_{\emptyset}[\chi_R] } \\
\tilde{g}_{m}(x_{m}) & := & \frac{ \mbox{M}_{m}[\chi_R . g](x_m) }{\mbox{M}_{m}[\chi_R](x_m)}- \tilde{g}_{\emptyset}  \nonumber \\
\tilde{g}_{mn}(x_m,x_n) & := & \frac{\mbox{M}_{mn}[\chi_R . g](x_m,x_n)}{\mbox{M}_{mn}[\chi_R](x_m,x_n)}   - \tilde{g}_m(x_m)- \tilde{g}_n(x_n) - \tilde{g}_{\emptyset}   \nonumber \\ \tilde{g}_{mm}(x_m,x_m) &=& 0 \nonumber .
\end{eqnarray}
We observe that this system is a generalization of (\ref{HDMR.approximation}) for an arbitrary approximation domain $R = \mbox{dom}(g) \subset X$. 

A new problem, however, arose as formulae (\ref{HDMR.approximation_mean}) are not valid any more in this general setting. As we have already indicated, these identities are beneficial in Section \ref{Section.HDMRminimization}, so we need to readjust all the components $\tilde{g}_\emptyset$, $\tilde{g}_m(x_m)$ and $\tilde{g}_{mn}(x_m,x_n)$ to satisfy (\ref{HDMR.approximation_mean}). Fortunately, this can be done easily without disturbing their optimality. We shift each component by the respective auxiliary constant $\sigma_\emptyset, \sigma_m, \sigma_{mn}$ in such a way that (\ref{HDMR.approximation_mean}) holds again. Formally, we define
\begin{eqnarray} \label{wHDMR.shifting}
\sigma_m & := & \sum_{x_m \in X_m} \tilde{g}_{m}(x_{m}) \\
\sigma_{mn} & := & \sum_{x_m \in X_m} \sum_{x_n \in X_n} \tilde{g}_{mn}(x_m,x_n)  \nonumber \\ \sigma_{mm} &:= & 0 \nonumber,
\end{eqnarray}
and then
\begin{eqnarray} \label{wHDMR.shifting0}
\sigma_\emptyset & := &  \sum_{m=1}^d \sigma_m + \frac{1}{2} \, \sum_{m,n=1}^d \sigma_{mn}.
\end{eqnarray}
Next, we redefine HDMR components as
\begin{eqnarray} \label{wHDMR.final}
\tilde{g}_{\emptyset} & := & \tilde{g}_{\emptyset} + \sigma_\emptyset \\
\tilde{g}_{m}(x_{m}) & := & \tilde{g}_{m}(x_{m})  - \sigma_m \nonumber \\
\tilde{g}_{mn}(x_m,x_n) & := & \tilde{g}_{mn}(x_m,x_n)  - \sigma_{mn}\nonumber .
\end{eqnarray}
The values of all $\sigma_m$ and $\sigma_{mn}$ determined by  (\ref{wHDMR.shifting}) now ensure the validity of (\ref{HDMR.approximation_mean}), and formula (\ref{wHDMR.shifting0}) guarantees that the overall shift of values of $\tilde{g}(\V{x})$ is nullified, see (\ref{GTILDE}), and so (\ref{wHDMR.final}) does not affect the optimality of $\tilde{g}(\V{x})$.

Even though equations (\ref{wHDMR.approximation}) and (\ref{wHDMR.final}) seem to be more complicated, their  computational complexity is similar to the full domain case (\ref{HDMR.approximation}). Therefore, we will refer to this more general result throughout this article. When $dom(g) = X$, both these approaches are equivalent. 
\section{Fast Minimization of HDMR}
\label{Section.HDMRminimization}
In this section, the main novelty of this article is developed. The key ingredient of the proposed approximate dynamic programming technique is a fast approximate minimization of functions in HDMR form.  We consider function $\tilde{g}(\V{x},\V{z})$, $\mbox{dom}(\tilde{g}) = X \times Z$, 
having the following structure
\begin{equation} \label{FormalHDMRMinimizationDEF}
\tilde{g}(\V{x},\V{z}) = \frac{1}{2} \sum_{m,n=1}^\mu  \tilde{g}_{mn}(z_m, z_n) +
\sum_{m=1}^\mu \tilde{g}_{m}(z_{m}) + \sum_{m=1}^\mu  \sum_{n=1}^\kappa \tilde{g}_{\mu+n,m}(x_n, z_m),
\end{equation}
where we denoted by $\kappa$ and $\mu$ the dimension of $X$ and $Z$,  respectively. This function corresponds to full HDMR of $\tilde{g}(\V{x},\V{z})$ without all HDMR components independent of $\V{z}$. Since we are interested in a point-wise minima of $\tilde{g}(\V{x},\V{z})$,
\begin{equation} \label{FormalHDMRMinimization}
p(\V{x}) := \MIN{\V{z} \in Z}{\tilde{g}(\V{x},\V{z})},
\end{equation}
the previous restriction on components of $\tilde{g}(\V{x},\V{z})$ is without loss of generality and it considerably eases the notation.

Regardless of a specific choice of $\V{x} \in X$, the parametrized minimization in (\ref{FormalHDMRMinimization}) is equivalent to the search for the clique of the minimal weight in a complete multi-partite edge-weighted graph \cite{Matousek1998_Invitation}. To show it, identify different $Z_m$ as partite sets of the graph, $z_m \in Z_m$ as vertices in particular partite set $Z_m$ and $\tilde{g}_{mn}(z_m, z_n) $ as weight of edge between vertices $z_m \in Z_m$ and $z_n \in Z_n$ taken from distinct partite sets with $\tilde{g}_{mm} = 0$, as we claimed in (\ref{HDMR.approximation_mean}). The additional weights of vertices $\tilde{g}_m(z_m)$ and $\tilde{g}_{\mu + n,m}(x_n,z_m)$, the latter parametrized by $x \in X$, can be simply added to the weights of proper edges.  This problem is known to be NP-hard \cite{Karp1972} and as it plays a role of repeatedly solved subproblem here, we search only for an approximate solution of (\ref{FormalHDMRMinimization}).  
\subsection{Problem Reformulation}
At the moment, it is fruitful to rewrite function $\tilde{g}(\V{x},\V{z})$ in a more convenient form. For a finite set $B$ and $i \in \{ 1, \ldots, | B| \}$ we denote $B[i]$ the $i$-th element of $B$. Then, for all $m,n \in \{ 1,\ldots, \mu \}$ we define matrices $\V{F^{mn}}$ in this way
\begin{equation} \label{defF}
F^{mn}_{ij} := \tilde{g}_{mn}(Z_m[i],Z_n[j]).
\end{equation}
In the same manner, we define matrices $\V{G^{mn}}$
\begin{equation} \label{defG}
G^{mn}_{ij} := \tilde{g}_{mn}(Z_m[i],X_{n}[j])
\end{equation}
for all $m \in \{ 1,\ldots, \mu \}$ and $n \in \{ 1,\ldots, \kappa \}$ and vectors $\V{h^m}$
\begin{equation} \label{defP}
h^m_i := \tilde{g}_{m}(Z_m[i])
\end{equation}
for all $m \in \{ 1,\ldots, \mu \}$.  Further, we compose all matrices $\V{F^{mn}}$ into one matrix $\V{F}$ with $\V{F^{mn}}$ being the $mn$-th subblock of $\V{F}$. Similarly, we create matrix $\V{G}$ out of matrices $\V{G^{mn}}$ and vector $\V{h}$ consisting of subvectors $\V{h^m}$. Thus, we obtain a concise reformulation of  $\tilde{g}(\V{x},\V{z})$ 
\begin{equation} \label{concisew}
\gamma(\V{u},\V{v}) := \frac{1}{2} \, \V{v}^T \V{F} \V{v}  + \V{h}^T \V{v} + \V{u}^T  \V{G} \V{v},
\end{equation}
where the only question left is to clarify the relation between  vectors $\V{u}$, $\V{v}$, and the original variables $\V{x} \in X$, $\V{z} \in Z$, respectively. 

We define
\begin{equation} \label{defDim}
\theta :=  \sum_{m=1}^\mu |Z_m|
\end{equation}
and follow the logic of the previous construction to deduce the structure of the newly introduced vector $\V{v} \in \mathbb{R}^\theta$. We see that it consists of $\mu$ subvectors 
\begin{equation} \label{z.equivalencea}
\V{v^{\,m}} \in \{0, 1 \}^{|Z_m|},
\end{equation}
which are related to coordinates $z_m \in Z_m$ of the original variable $\V{z} \in Z$ as
\begin{equation} \label{z.equivalence} 
v^m_i := 1 \Longleftrightarrow z_m = Z_m[i], \qquad v^m_i := 0 \quad \mbox{otherwise}.
\end{equation}
The relation of vector $\V{u}$ to the original parameter $\V{x} \in X$ is analogous. Such constructions of $\V{v}(\V{z})$ and $\V{u}(\V{x})$ guarantee that 
\begin{equation} \label{gamma.g.equivalence} 
\gamma(\V{u}(\V{x}),\V{v}(\V{z})) = \tilde{g}(\V{x},\V{z}),
\end{equation}
for all $(\V{x},\V{z}) \in X \times Z$, and thus the evaluation of $p(\V{x})$, see (\ref{FormalHDMRMinimization}), is fully equivalent to minimization of $\gamma(\V{u}(\V{x}),\V{v})$ with respect to all vectors $\V{v}$ obeying (\ref{z.equivalence}). Therefore, the latter problem is also a NP-hard problem. It is, however, more amenable to the relaxation technique developed further.
\subsection{Trust Region Based Relaxation} \label{TRBasedMinimisation}
We observe that each $\V{x} \in X$ in (\ref{FormalHDMRMinimization}) yields a different value of parameter $\V{u}$ in (\ref{concisew}) while matrix $\V{F}$  remains unchanged. Thus, we can afford some intensive matrix  preprocessing in order to exploit the repetitive nature of this minimization.
That is why we turn our attention to the trust region problem \cite{Sorensen1982_Trust_Region} which permits fast exact solution even for an indefinite matrix $\V{F}$.  To match the form of the trust region problem, we have to relax constraints (\ref{z.equivalence}) into $\| \V{v}\| = r$ with $r>0$ specified lately. Thus, we obtain problem
\begin{equation} \label{TR_forW}
\MIN{\| v\| = r} \left \{ \frac{1}{2} \, \V{v}^T \V{F} \V{v} + \V{h}^T  \V{v} +  \V{u}^T \V{G}  \V{v}  \right \}.
\end{equation}
The only question left is to adjust the diameter $r$ properly. 

We can set $r^2=\mu$ immediately, as each feasible vector $\V{v}$ of the original problem consists of $\mu$ subvectors $\V{v^m}$ of unit norm, see (\ref{z.equivalence}). Yet there is a possibility of obtaining a tighter relaxation.  By the definition of matrices $\V{F}$, $\V{G}$ and vector $\V{h}$, see (\ref{defF}), (\ref{defG}) and (\ref{defP}), respectively, and by zero mean of all HDMR components derived in (\ref{HDMR.approximation_mean}) and (\ref{wHDMR.final}), we observe that the minimized criteria in (\ref{TR_forW}) do not depend on the average value of any subvector $\V{v^{m}}$ of $\V{v}$. Hence, we may shift all elements of each $\V{v^{m}}$ by a constant factor $- \frac{1}{|Z_m|}$ and thus rewrite constraint (\ref{z.equivalence}) as
\begin{equation} \label{z.equivalenceee}
v^m_i := 1 - \frac{1}{|Z_m|} \Longleftrightarrow z_m = Z_m[i], \qquad v^m_i := - \frac{1}{|Z_m|} \quad \mbox{otherwise},
\end{equation}
and the value of $\gamma(\V{u},\V{v})$ remains unchanged. This observation suggests adjusting a slightly smaller diameter $r$ in this manner
\begin{equation} \label{choiceR}
r^2 := \sum_{m =1}^\mu \left \{ \left ( 1 - \frac{1}{|Z_m|} \right )^2 + \sum_{i=2}^{|Z_m|} \frac{1}{|Z_m|^2} \right \} = \mu - \sum_{m =1}^\mu \frac{1}{|Z_m|},
\end{equation}
which corresponds to the norm of any feasible solution satisfying constraint (\ref{z.equivalenceee}). Thus, we obtained as tight relaxation of the original problem as possible and we are ready to solve the trust region problem (\ref{TR_forW}).

From a wide spectra of solution methods of the trust region problem, see \cite{RojasSantos2000_MatrixFreeTRP}, and references therein, we choose one which is computationally expensive for a one step minimization, but effective in our repetitive setting. At first, we find ortoghonal matrix $\V{U}$ such that 
\begin{equation} \label{def_UD}
\V{F} = \V{U}^T \V{D} \V{U}
\end{equation}
holds with diagonal matrix $\V{D}$ having its diagonal composed of all eigenvalues ordered from the lowest one to the highest one. Then, for a particular $\V{u}$ we define 
\begin{equation} \label{def_b}
\V{b} := \V{U} \V{h} +  \V{U} \V{G}^T \V{u}, 
\end{equation}
and we find solution $\hat{\V{v}}$ of (\ref{TR_forW}) according to
\begin{equation} \label{exact_minimizer_w}
\hat{\V{v}} := - \V{U}^T ( \V{D} - \lambda \, \mathbb{I})^{-1} \V{b},
\end{equation}
where $\mathbb{I}$ is unit matrix and $\lambda \in (-\infty, D_{kk} )$ solves one-dimensional equation 
\begin{equation} \label{calculating.lambda}
\sum_{i=k}^\theta \left ( \frac{b_i}{D_{ii} - \lambda} \right )^2 = r^2,
\end{equation}
with an index of the first non-zero element of $\V{b}$ denoted by $k \in \{1,\ldots,\theta\}$. Then, precisely one such $\lambda$ exists
and can, for instance, be computed by the Newton's method.  
A more detailed discussion is to be found in 
 \cite{RojasSantos2000_MatrixFreeTRP,Olsson2007_Large_Scale_BIQ,Sorensen1982_Trust_Region,Busygin}. 

In some practical problems, matrix $\V{F}$ in (\ref{TR_forW}) may be zero or may have a very small norm. Then, we may either use some different approach, e.g. linear integer programming \cite{Schrijver1998_LIP}, or we may solve (\ref{TR_forW}) analytically with the optimal choice of $\hat{v}$ determined by formula
\begin{equation} \label{nullfvariant}
\hat{\V{v}} = - \frac{ \| r \| }{\| \V{h} +  \V{G}^T \V{u} \|} \left ( \V{h} +  \V{G}^T \V{u}  \right ).
\end{equation}
\subsection{Estimate of the Exact Minimizer}
\label{estimate_of_minimizer}
At the moment, we briefly summarize the previous procedure. We related $\V{v}(\V{z}) \in \mathbb{R}^\mu$ to each $\V{z} \in Z$ by (\ref{z.equivalence}), and also $\V{u}(\V{x}) \in \mathbb{R}^\kappa$ to $\V{x} \in X$ in a similar manner. Next, we found the exact minimizer $\hat{\V{v}} \in \mathbb{R}^\theta$ of the relaxed problem (\ref{TR_forW}), which is in fact parametrized by $\V{x} \in X$ as $\hat{\V{v}} = \hat{\V{v}}(\V{u}(\V{x})) = \hat{\V{v}}(\V{x})$. Such $\hat{\V{v}}(\V{x})$ generally does not correspond to any feasible solution $\V{z} \in Z$ of the original problem (\ref{FormalHDMRMinimization}). Yet, we may still use the knowledge of $\hat{\V{v}}(\V{x})$  to estimate the value of $p(\V{x})$. 

First, we easily obtain a lower bound 
\begin{equation} \label{lowerp}
\underline{p}(\V{x}) := \gamma(\V{u}(\V{x}),\hat{\V{v}}(\V{x})).
\end{equation}
Indeed, if we compare the derivation of (\ref{TR_forW}) with the original problem (\ref{FormalHDMRMinimization}), we realize that $\gamma(\V{u}(\V{x}),\hat{\V{v}}(\V{x}))$ minimizes the same criteria with respect to a larger set. Therefore, we have $\underline{p}(\V{x}) \le p(\V{x})$ for all $\V{x} \in X$. This lower bound $\underline{p}(\V{x})$ is, however, problematic. It gives only poor estimates on $p(\V{x})$ as we show in a numerical experiment in Section \ref{MinRandFun}. 

On that account, we develop a more accurate upper estimate on $p(\V{x})$ now. We simply interpret each value $\hat{v}_i^{m}(\V{x})$ as an indicator of subobtimality of the related element $Z_m[i]  \in Z_m$. In other words, the higher the element $\hat{v}_i^{m}(\V{x})$ is, the lower cirteria $\tilde{g}(\V{x},z_1,\ldots,z_\mu)$ we may expect when adjusting $z_m$ to $Z_m[i]$. One can came up with many different ways of such \QT{rounding} of  $\hat{\V{v}}(\V{x})$ to some $\V{z} \in Z$, and thus there is not any guarantee that the following heuristic is the best possible. 

From now on, we again omit parameter $\V{x} \in X$ in the notation for the sake of simplicity. We start with normalizing vector  $\hat{\V{v}}$ in two setps. We shift it to be non-negative 
\begin{equation}
\hat{\V{v}}:= \hat{\V{v}} - \MIN{i  \in \{1,\ldots,\theta\}}{\hat{v}_i}, 
\end{equation}
and then we rescale all its subvectors $\V{\hat{v}^m}$, $m \in \{ 1,\ldots,\mu\}$, as follows
\begin{equation} \label{maximizers}
\V{\hat{v}^m} :=  \V{\hat{v}^m} \,\, /  \MAX{i \in \{1,\ldots,|Z_m|\}}{\hat{v}^{m}_i}.
\end{equation}
Thus, for all $m$ there is at least one element of $\V{\hat{v}^m}$ equal to $1$, and for all $i \in \{1,\ldots,|Z_m|\}$ it holds that $\hat{v}^m_i \in [0,1]$. 
Further, we define function $q(\V{z})$ indicating the quality of a particular $\V{z} \in Z$ (with respect to an implicit parameter $\V{x} \in X$)
\begin{equation}  \label{w_z_def}
q(\V{z}) := \prod_{m=1}^\mu \hat{v}^m_{i_m} \quad \mbox{where}\quad  \V{z} =\left (Z_1[i_1],\ldots,Z_{\kappa}[i_{\kappa}] \right ).
\end{equation}
From non-negativity of $\hat{\V{v}}$ we observe that $q(\V{z}) \in [0,1]$, and the maximum of $q(\V{z})$ with respect to $\V{z} \in Z$ is equal to $1$ by (\ref{maximizers}).  Then, we define
\begin{equation} \label{defZm}
Z^\phi := \left \{ \V{z} \in Z :  q(\V{z}) \ge \phi \right \}  \end{equation}
for any $\phi \in [0,1]$. Thus, $Z^0 = Z$ and $Z^1$ contains only such $\V{z} \in Z$ that all the corresponding $\hat{v}^m_{i_m} $ are maximizers used in the denominator in (\ref{maximizers}).  We note that $Z^\phi$ can be enumerated in a component-wise manner using (\ref{w_z_def}) without passing the whole $Z$. 
Then, we substitute $Z^\phi \subset Z$ for $Z$ in  (\ref{FormalHDMRMinimization}), and we find an upper bound $\overline{p}^\phi(\V{x})$ on $p(\V{x})$
\begin{equation} \label{def_overlined_p}
\overline{p}^\phi(\V{x}) := \MIN{\V{z} \in Z^\phi}{\tilde{g}(\V{x},\V{z})},
\end{equation}
by enumerating $\tilde{g}(\V{x},\V{z})$ for all $\V{z} \in Z^\phi$. The lower the value of $\phi$ we choose, the larger the $Z^\phi$ that we obtain and the tighter the upper bound $\overline{p}^\phi(\V{x})$ we find; nonetheless, at the price of slower enumeration in (\ref{def_overlined_p}). 

Once the diagonalization in (\ref{def_UD}) is done, it is in fact easy to compute $\overline{p}^\phi(\V{x})$ for any $\V{x} \in X$. We construct $\V{u}(\V{x})$ by the one-to-one correspondence (\ref{z.equivalence}), then we compute vector $\V{b}(\V{x})$ according to (\ref{def_b}), find the related value of $\lambda(\V{x})$ following (\ref{calculating.lambda}), and finally calculate candidate $\hat{\V{v}}(\V{x})$ which enters the already introduced procedure that leads to $\overline{p}^\phi(\V{x})$ defined by (\ref{def_overlined_p}).  Thus, we found  approximate minima of a general function $\tilde{g}(\V{x},\V{z})$ in HDMR form over $\V{z} \in Z$ for all parameters $\V{x} \in X$. This permits us to apply HDMR to effectively approximate the Bellman equation  in Section \ref{Sec.ADPHDMR}.
\subsection{Minimization of a Random Function} 
\label{MinRandFun}
\begin{figure}
\includegraphics[scale=0.82]{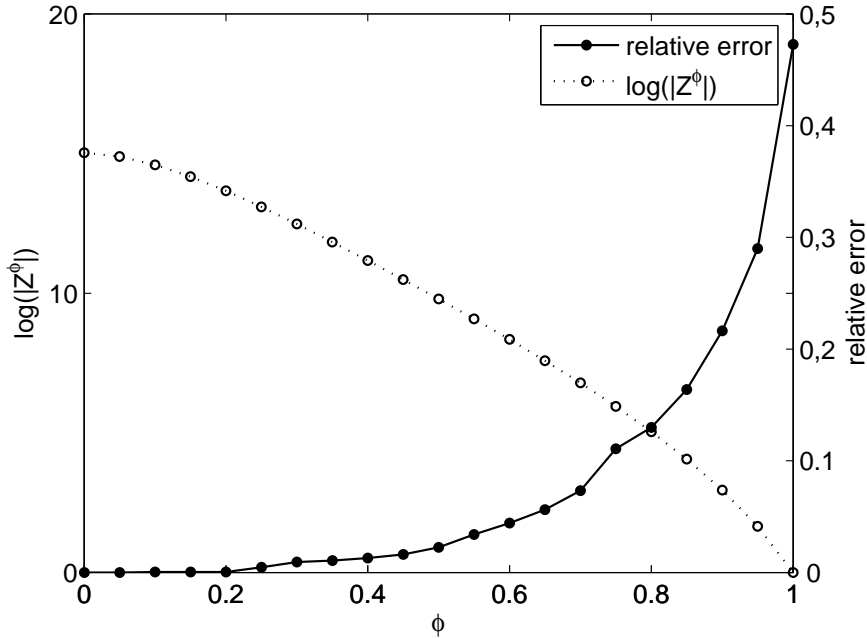}
\caption{Value of $\mbox{log}(|Z^\phi|)$ and a relative error of $\overline{p}^\phi$ plotted against various values of $\phi$. The relative error is the distance of $\overline{p}^\phi$ from the minimum of $\tilde{g}(z)$ rescaled and shifted in such a way that exact minimum corresponds to $0$ whereas the average value of the minimized criteria correspods to $1$. The depicted results were averaged over 20 different realizations of matrix $\V{F}$ and vector $\V{h}$} 
\label{Fig1}
\end{figure}

Now, we dedicate a short section to a numerical verification of the previously introduced technique. We solved problem (\ref{FormalHDMRMinimization}) exactly for a random function $\tilde{g}(\V{x},\V{z})$. For the sake of simplicity, we omitted parameter $\V{x}$ and set $\V{G} = 0$ in (\ref{TR_forW}). Next, we choose the minimization domain $Z = \{1,\ldots,150\}^3$, we generated HDMR components $\tilde{g}_{mn}$ randomly with values chosen from uniform distribution on interval $[0,1]$ and finally we adjusted them to satisfy (\ref{HDMR.approximation_mean}). Then, we found a lower estimate on minima $\underline{p}$ according to (\ref{lowerp}) and upper estimates on minima $\overline{p}^\phi$ for various choices of parameter $\phi$ following (\ref{def_overlined_p}). All results were averaged with respect to 20 random samples of $\V{F}$ and $\V{h}$ and depicted in Fig. \ref{Fig1}. The relative error of upper bound $\overline{p}^\phi$ is defined as the distance from minimum of $\tilde{g}(\V{z})$ rescaled and shifted in such a way that the exact minimum corresponds to $0$ whereas the average value of the minimized criteria corresponds to $1$. 
We observe that the lower the value of $\phi$ is, the better the approximation we obtain as we expected. On the other hand, there was a linear grow of $\mbox{log}(|Z^\phi|)$ when decreasing $\phi$. We suppose that a detailed elaboration of this relation could serve as a basis for an error estimation heuristics. Concerning the lower bound, we obtained $\underline{p} = -5.55$ holding the same scale as previously,
whereas the worst upper bound $\overline{p}^1 = 0.47$ is almost $12$ times closer to the exact minimum $0$.
As both have similar computational complexity, we omit lower bound estimate $\underline{p}$ from further considerations.

These experiments were carried out on CPU Intel Core i3, 2.10 GHz with 4GB of RAM in Matlab 7.
It took $169$ seconds to find the exact minimum, whereas the average time necessary to diagonalize matrix $\V{F}$ was $1.3$ seconds. We note that this matrix diagonalization is done only once in the full setting of (\ref{FormalHDMRMinimization}), whereas the time necessary for exact minimization of $\tilde{g}(\V{x},\V{z})$ for each $\V{x} \in X$ is still the same.
\section{Approximate DP based on HDMR}
\label{Sec.ADPHDMR}
This is the right time to briefly introduce the decison-making theory. A~decision-making task stands for selecting a decision-maker's strategy
in order to reach his aim with respect to the part
of the world (so-called system). The decision maker observes or influences the
system over a finite decision making horizon $\tau < \infty$. Value $\V{y}_{t} \in \V{y}_t$, $t \in T = \{ 1, \ldots, \tau \}$, provides the decision maker with all the knowledge influencing the future behaviour of the system. Thus, $\V{y}_t$ includes the current state of the system together with other external data observed up to time instant $t$. Nonetheless, we will reference $\V{y}_t$ simply as a state of the system. Next, the decisions (actions) of a decision-maker are denoted as
$\V{a}_{t} \in A_t$. A strategy is a collection of mappings of the current state $\V{y}_{t-1} \in Y_{t-1}$ into the choice of the next
decision $\V{a}_t \in A_t$; for the optimal strategy we use symbols $\{ \hat{\V{a}}_{t}(\V{y}_{t-1}) \}_{t \in T}$.
To formalize the decision-maker's aims, a concept of the additive loss function is used, $l_t (\V{a}_{t},\V{y}_{t})$, depending on the current action $\V{a}_t$ and system state $\V{y}_t$. The involved system is described
in a probabilistic manner by the following collection of pdfs called the outer Markov model of a system
\begin{equation} 
\label{SystemModel}
\left \{ f_t(\V{y}_t|\V{a}_{t},\V{y}_{t-1}) \right \}_{t \in T}. \end{equation} 

For the expected value of variable $x$ conditioned by $y$ we use
\begin{equation} 
\E{x}{y} := \int_X x\,f(x|y)\,\mbox{d}x.
\end{equation}
Knowing the collection of loss functions $\{ l_t (\V{a}_{t},\V{y}_{t}) \}_{t=1}^\tau$ together with the system model (\ref{SystemModel}), the optimal strategy $\{ \hat{\V{a}}_{t}(\V{y}_{t-1}) \}_{t \in T}$ is fully determined by the Bellman function
\begin{equation} 
V_{t-1}(\V{y}_{t-1}) = \MIN{\V{a}_t \in A_t} \E{l_t(\V{a}_{t},\V{y}_{t})+V_{t}(\V{y}_{t})}{\V{a}_{t},\V{y}_{t-1}},
 \label{Bellman} 
\end{equation} 
which has to be recursively evaluated at all times $t \in T$ with the boundary condition $V_{\tau} = 0$. As this standard form of the Bellman equation (\ref{Bellman}) is not convenient to our purposes, we rewrite it in an equivalent form
\begin{eqnarray} \label{RewrittenBellman}  
E_{t}(\V{a}_t,\V{y}_{t-1}) & = & \E{l_t(\V{y}_t,\V{a}_t)+ \MIN{\V{a}_{t+1} \in A_{t+1}}
 E_{t+1}(\V{a}_{t+1},\V{y}_{t})}{\V{a}_t,\V{y}_{t-1}} \\
E_{\tau+1} & = & 0. \nonumber
\end{eqnarray} 
Then, $E_{t+1}(\V{a}_{t+1},\V{y}_t)$ is the expected loss-to-go provided we choose action $\V{a}_{t+1}$ in the system state $\V{y}_{t}$. In this setting, the optimal strategy $\hat{\V{a}}_{t}(\V{y}_{t-1})$ is composed of actions satisfying 
\begin{equation} \label{optimal_choice}
\hat{\V{a}}_{t}(\V{y}_{t-1}) := \ARGMIN{\V{a}_t \in A_t} E_t(\V{a}_t,\V{y}_{t-1}).
\end{equation} 
\subsection{Offline Part - Approximate Evaluation  of $E_t$} 
\label{Offline}
Now, we are prepared to apply both HDMR developed in Section \ref{Section.HDMR} and fast approximate minimization of functions in HDMR form, see Section \ref{Section.HDMRminimization}, to effectively approximate $E_{t}(\V{a}_t,\V{y}_{t-1})$ defined by (\ref{RewrittenBellman}). This part of algorithm is the most demanding concerning the computational complexity. Thus, function $E_{t}(\V{a}_t,\V{y}_{t-1})$ is typically computed offline, stored as a look-up table (in our case in HDMR form), and then used during the online part of a decision-making algorithm to find the approximated optimal action by using $(\ref{optimal_choice})$. The proposed algorithm runs in the backward manner analogously to the evaluation of the exact Bellman equation (\ref{Bellman}). 

We denote the approximated loss-to-go function by $\tilde{E}_t$ even though for $t < \tau $ it is not the exact HDMR of $E_t$ . For the first step, $t = \tau$, we rewrite (\ref{RewrittenBellman}) as
\begin{equation} \label{horizon_e}
E_{\tau}(\V{a}_\tau,\V{y}_{\tau-1}) = \E{l_{\tau}(\V{y}_\tau,\V{a}_\tau)}{\V{a}_\tau,\V{y}_{\tau-1}}.
\end{equation} 
To obtain all HDMR components $\tilde{E}_{\tau,\emptyset},\tilde{E}_{\tau,m},\tilde{E}_{\tau,mn}$ of  $E_{\tau}(\V{a}_\tau,\V{y}_{\tau-1})$, we evaluate $E_{\tau}(\V{a}_\tau,\V{y}_{\tau-1})$ for each pair $(\V{a}_\tau,\V{y}_{\tau-1}) \in A_\tau \times Y_{\tau-1}$ and add the resulting value to proper sums in (\ref{wHDMR.final}). 

Next, suppose we know all $\tilde{E}_{t+1,\emptyset}$, $\tilde{E}_{t+1,m}$, $\tilde{E}_{t+1,mn}$ and we want to find an approximation of $E_t$ in the form of HDMR. Substituting $\tilde{E}_{t+1}$ into (\ref{RewrittenBellman}) we have
\begin{equation} \label{recursive_approx}
E_{t}(\V{a}_t,\V{y}_{t-1}) \approx \E{l_t(\V{y}_t,\V{a}_t)+ \MIN{\V{a}_{t+1} \in A_{t+1}}
\tilde{E}_{t+1}(\V{a}_{t+1},\V{y}_{t}))}{\V{a}_t,\V{y}_{t-1}}.
\end{equation} 
This suggests defining $\tilde{E}_{t}$ as HDMR of the expression on the right-hand side, or at least as HDMR of some approximation of this expression. On that account we denote 
\begin{equation} \label{one_step_minimization}
\pi_t(\V{y}_t) := \MIN{\V{a}_{t+1} \in A_{t+1}}
\tilde{E}_{t+1}(\V{a}_{t+1},\V{y}_{t})
\end{equation} 
and search for its upper bound $\overline{\pi}_t^\phi(\V{y}_t)$ following the instructions of Section \ref{Section.HDMRminimization}. The choice of an auxiliary parameter $\phi \in [0,1]$  determining the precision of the upper bound estimate is discussed at the end of this section.
Looking at (\ref{FormalHDMRMinimization}), we identify $\tilde{g} = \tilde{E}_{t+1}$, $X = \V{y}_t$ and $Z = A_{t+1}$.  We note that all the HDMR components of  $\tilde{E}_{t+1}$ that depend only on $\V{y}_{t}$ may be directly interchanged with minimization in (\ref{one_step_minimization}) and thus not considered at the moment. Based on the knowledge of such $\tilde{E}_{t+1,\emptyset}$, $\tilde{E}_{t+1,m}$ and $\tilde{E}_{t+1,mn}$ that depend on $\V{a}_{t+1}$, we construct matrices $\V{F}_t$, $\V{G}_t$ and vector $\V{h}_t$ according to (\ref{defF}), (\ref{defG}) and (\ref{defP}), and we formulate the relaxed problem (\ref{TR_forW}). Then, we find its exact minimizer $\hat{\V{v}}_t(\V{y}_t)$ in a direct analogy to (\ref{exact_minimizer_w}) with matrix diagonalization
\begin{equation} \label{diagonalization_bellman}
\V{F}_t = \V{U}_t^T \V{D}_t \V{U}_t
\end{equation} 
involved. The diagonalized matrix $\V{F}_t$ is typically small and does not grow much with $t$ as its size (\ref{defDim}) corresponds to the space of actions $\V{a}_t$. Knowing $\hat{\V{v}}_t(\V{y}_t)$, we calculate an upper bound on minimum applying procedure (\ref{def_overlined_p}), and finally we add (restore) all HDMR components of $\tilde{E}_{t+1}$ depending only on $\V{y}_{t}$. Thus, we obtained an upper bound on minimum of $\pi_t(\V{y}_t)$.
We note that diagonalization (\ref{diagonalization_bellman}) is carried out just once for each time step $t$, and so we can effectively evaluate $\overline{\pi}_t^\phi(\V{y}_t)$ for all $\V{y}_t \in \V{y}_t$. Now, we find $\tilde{E}_t(\V{a}_t,\V{y}_{t-1})$ by evaluating 
the right-hand side of the following formula 
\begin{equation} \label{aprox_final}
\tilde{E}_{t}(\V{a}_t,\V{y}_{t-1}) \approx \E{l_t(\V{y}_t,\V{a}_t)+ \overline{\pi}^\phi_t(\V{y}_t)}{\V{a}_t,\V{y}_{t-1}}
\end{equation} 
for each pair $(\V{a}_t,\V{y}_{t-1}) \in A_t \times Y_{t-1}$ and add the resulting value to proper sums in (\ref{wHDMR.final}) immediately. Thus, we construct  all HDMR components $\tilde{E}_{t,\emptyset}$, $\tilde{E}_{t,m}$ and $\tilde{E}_{t,mn}$, avoiding the full dimensional representation of $\tilde{E}_t$. 

Finally, we repeat the whole procedure to recursively compute function $\tilde{E}_t(\V{a}_t,\V{y}_{t-1})$ for all $t \in T$. Once the calculation of each particular $\tilde{E}_t$ is finished, we can completely remove all components of $\tilde{E}_{t+1}$ independent of $a_{t+1}$ non-affecting the suboptimal strategy computed in the next section.
\subsection{Online Part - Approximate Minimization of $\tilde{E}_t$} 
\label{Online}
The previously described part of the algorithm has to be implemented in advance, or \QT{off-line} manner because of high computational demands. As functions $\{ \tilde{E}_t(\V{a}_t, \V{y}_{t-1}) \}_{t \in T}$ are stored only in the form of HDMR, it is possible to take larger decision horizons $\tau$ into consideration. Nonetheless, we still have to choose an approximated (suboptimal) action $\tilde{a}_t$ in the real time, or \QT{on-line} manner. Then, the previously observed system state $\V{y}_{t-1}$ is fixed and so we solve just one minimization problem in each time step $t$ in opposite to the recursive evaluation of (\ref{recursive_approx}). Substituting $\tilde{E}_t$ into (\ref{optimal_choice}), we define 
\begin{equation} \label{suboptimal_choice}
\tilde{\V{a}}_t(\V{y}_{t-1}) := \ARGMIN{\V{a}_t \in A_t} \tilde{E}_t(\V{a}_t,\V{y}_{t-1}).
\end{equation} 
We note that $\tilde{\V{a}}_t(\V{y}_{t-1})$ does not stand for HDMR approximation of $\V{\hat{a}}_t(\V{y}_{t-1})$ defined by (\ref{optimal_choice}). 

There are many ways how to find $\tilde{\V{a}}_t$, or at least some its approximation. An interesting choice can be a trust region based relaxation as we may exploit our previous calculations. We may represent HDMR components of $\tilde{E}_t(\V{a}_t,\V{y}_{t-1})$ in the basis obtained 
in (\ref{diagonalization_bellman}). If we store all matrices $\V{U_t}$, $\V{D_t}$, and also matrices $\V{G_t}$ and vectors $\V{h}_t$  involved in approximate minimization of $\pi_t(\V{y}_t)$ defined by (\ref{one_step_minimization}), we may find approximate minimizer of (\ref{suboptimal_choice}) in accordance with Section \ref{Section.HDMRminimization} again. However, even some more accurate technique may be used in one-shot only  minimization (\ref{suboptimal_choice}). Any algorithm for binary quadratic programming \cite{Olsson2007_Large_Scale_BIQ} may be applied to solve (\ref{suboptimal_choice}) via equivalent reformulation (\ref{concisew}) constrained by (\ref{z.equivalence}). For smaller sets $A_t$, we can find even exact value of $\V{\tilde{a}}_t \in A_t$ by direct enumeration of (\ref{suboptimal_choice}). We decided to use this most accurate approach in Section \ref{NBandit} in order to show the extent to which $\tilde{E}_t(\V{a}_t,\V{y}_{t-1})$ in the form of HDMR may be compared with exact value of $E_t(\V{a}_t,\V{y}_{t-1})$.

\subsection{N-armed Bandit Problem} 
\label{NBandit}
As an ilustrative example, we propose here an approximate solution to the $N$-armed bandit problem, which was extremely important in approximate dynamic programming, see for instance \cite{Sutton1998_Reinforcement_Learning,Powell2007_ApproximateDP} and references therein. We compare its exact solution with HDMR based approximation. 

Conceive a game where the player has to choose between different options, e.g. levers of $N$-armed bandit, with numerical rewards chosen from various stationary probability distributions. The payoff probabilities of levers are fixed, yet unknown, and thus the player has to estimate them. Then the problem is to identify the most winning lever. Even though this problem could be formulated easily, it is a real issue for a longer game horizon as it is hard to balance exploration and exploitation. Winning in the first round does not imply that the player should stick to the same lever as it prevents learning of the payoff probability of other levers.

We considered game with $9$-armed bandit and decision making horizon of $\tau = 8$ steps to be able to compare approximated suboptimal strategies with the exact optimal strategy. Using the previous notation, $y_t$ stands for the observed value (payoff) $y_t \in Y = \{0,1\}$ and $\V{a}_t$ denotes the 
decision of a player in each time step $t \in T = \{1,\ldots,\tau\}$. The arms of the bandit are represented by two-dimensional space of actions, $\V{a}_{t} \in A = \{1,2,3\}^2$. The loss function 
\begin{equation} \label{losssf}
l_t(y_t,\V{a}_t) = -y_t
\end{equation} 
represents the aim of maximizing the payoff $y_{t}$ in each round of the game. Next, we introduce a sufficient statistic $\V{s}_t$, $\mbox{dom}(\V{s}_t)= Y \times A$, which compresses the previous game results in a small vector
\begin{equation} 
\V{s}_t(y,\V{a}) := \V{s}_{t-1}(y,\V{a}) + \delta_{y_t,y} \, \delta_{\V{a}_t,\V{a}},
\end{equation} 
with $\delta$ standing for standard Kronecker's symbol. Thus, $\V{s}_t(y,\V{a})$ counts how many times we observed a value $y$ after selecting an action $\V{a}$ in first $t$ rounds of the game. We set $\V{s}_0 = \V{0}$ for the moment. In fact, $\V{s}_t$ may be included into the system state $y_t$, but for the sake of simplicity we treat it separately here.  To compute the expected loss in (\ref{RewrittenBellman}), the knowledge of the Markov system model (\ref{SystemModel}) is necessary 
\begin{eqnarray}
f_t(y_t | \V{a}_t , \V{s}_{t-1} ) &=& \frac{\V{s}_{t-1}(y_t,\V{a}_t)+1}{\V{s}_{t-1}(y_t,\V{a}_t)+\V{s}_{t-1}(1-y_t,\V{a}_t)+2}.
\end{eqnarray}
This model was obtained using the technique of Bayesian estimation \cite{Peterka1981}. In the following experiment, the $9$-armed bandit was simulated using pseudo-random generator with fixed payoff probability matrix $\V{P}$ defined for  $\V{a} \in A$ as follows
\begin{equation} 
P_{ij} := \mbox{Prob}(y = 1 |\V{a} = [i,j] ).
\end{equation} 

During the experiment, it turned out that high-symmetry of $N$-armed bandit is unsuitable for our purposes. If the underlying payoff probability $\V{P}$ is completely unknown, and for the prior information it holds $s_0(y,\V{a}) = 0$ for all $y \in Y, \V{a} \in A$, then all the bandit arms have the same expected loss when averaged over all the possible system trajectories. Thus, $\V{F}_t$ corresponding to differences of the expected loss among various arms is equal to zero. We may still use the previously introduced algorithm, see the note near (\ref{nullfvariant}), but we would miss its most interesting part, i.e. the trust region based approximate minimization described in Section  \ref{TRBasedMinimisation}. 
We note that this high level of symmetry is very unlikely for a real-world problem.

Thus, we decided to slightly perturb the experiment to suppress its symmetry. We put a prior information on one arm, $s_0(0,[1,1]) = 1$, and in this setting we computed the exact values of $\{ E_t \}_{t \in T}$ following (\ref{RewrittenBellman}) and also all HDMR functions  $\{ \tilde{E}_t^\phi \}_{t \in T}$ according to Section \ref{Offline}. This time we explicitly stated that $\tilde{E}_t^\phi$ depends also on the value of $\phi$, see (\ref{aprox_final}). The disk space necessary to save $\{E_t\}_{t \in T}$ and each $\{ \tilde{E}_t^\phi \}_{t \in T}$ in Matlab .mat file was $2.3$~MB and $0.1$~MB, respectively. The optimal strategy was derived from $E_t$ using (\ref{optimal_choice}), and suboptimal strategies parametrized by $\phi$ were derived according to (\ref{suboptimal_choice}). 

\begin{figure}
\includegraphics[scale=0.82]{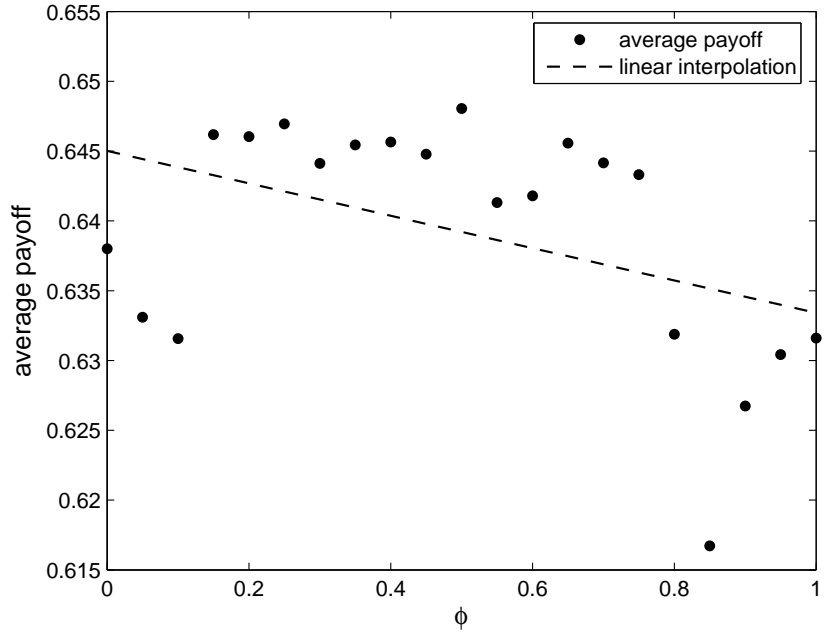}
\caption{Average payoffs obtained from
approximated strategies derived from $\tilde{E}_{t}^\phi$ for various $\phi \in [0,1]$. The average payoff of the exact optimal strategy derived from $E_t$ was $0.653$. 
These results are based on $20000$ simulated plays with $9$-armed bandit, each of them consisting of $\tau = 8$ steps. The payoff probabilities $P_{ij}$ of the bandit were chosen randomly from uniform distribution on interval $[0,1]$. The only exception was payoff probability $P_{11} = 0.1$, which was kept fixed to avoid complete symmetry of the problem as discussed in Section \ref{NBandit}}
\label{Fig2}
\end{figure}

All these strategies were used to simulate $20000$ plays with $9$-armed bandit, each of them consisting of $\tau = 8$ steps. The payoff probabilities $P_{ij}$ of the bandit were chosen randomly from uniform distribution on interval $[0,1]$ with the only exception of fixed payoff probability $P_{11} = 0.1$ corresponding to the only non-zero prior $s_0(0,[1,1])$. The average payoff of the optimal strategy was $0.653$, and the average payoffs obtained for various values of $\phi$ are depicted in Fig. \ref{Fig2}. The strategy derived from $E_t^1$ was rather sucessful, it gained $0.632$ on average. It indicates the practical applicability of the less acccurate approximation of $E_t$, when $\phi = 1$ and $Z^1$ contains typically just one element. Then, the whole estimating of the exact minimizer, see Section  
\ref{estimate_of_minimizer},  amounts only to \QT{rounding} of trust region problem minimizer to an approximate minimizer of HDMR. The precision of HDMR approximation itself may be deduced from the average payoff $0.638$ obtained for $\phi = 0$, which corresponds to the exact minimization in (\ref{one_step_minimization}). The closer the $\phi$ is to $0$, the closer $E_t^\phi$ is to $E_t$ by its definition. However, this monotonicity does not hold for the derived strategies. Yet, on average it holds again, see the interpolated line in Fig. \ref{Fig2}. The slope of this line is rather small; it means that in this particular problem the average payoff just slightly increases when decreasing $\phi$. It is in contrast with Fig. \ref{Fig1} where the upper bound estimate depended strongly on the minimization precision tuned by parameter $\phi$. Nonetheless, if we find upper bound $\overline{\pi}_{\tau-1}^\phi(y_{\tau-1})$ on (\ref{one_step_minimization}) for various $\phi$ and compare it with exact minimizer $\pi_{\tau-1}(y_{\tau-1})$, we obtain dependence on $\phi$ similar to that depicted in Fig.~\ref{Fig1}.
Thus, we observed better performance of strategies derived with $\phi$ close to $1$ than we can expect from the quality of upper bound estimates on $E_t^\phi$. This may be explained by some sort of systematic error produced by approximate minimization. Consider some fixed $\phi$. If all values of $E_t^\phi$ overestimate (or underestimate) values of $E_t$ by the same number, the approximate minimization would give inaccurate results, but both approximate and optimal strategies derived from $E_t^\phi$ and $E_t$, respectively, would be the same. However, more work has to be done to fully verify this conjecture, which is likely to be problem-dependent. 
\section{Conclusion}
\label{Section.Conclusion}
The aim of this work was to cope with both computational and  memory demands necessary to find and represent the optimal decision making strategy. The proposed variant of approximate dynamic programming based on HDMR is appealing for two reasons.  At first, this approximation considerably reduces memory demands, but, more importantly, it also enables a fast approximate minimization of the approximated Bellman function. Results of numerical simulation proved that the proposed variant of dynamic approximate programming is a viable technique. 

As for all the approximate methods surveyed at the beginning of Section \ref{Section.Introduction}, the one proposed in this article cannot be assigned to any of these classes directly. It is based on the Bellman function approximation; however, looking at its internal structure it may be considered also as an aggregation method where each HDMR component aggregates a different coordinates. Next, the point-wise construction of HDMR resembles the learning phase of the artifical neural networks, yet it is more straightforward. 

A bottleneck of the proposed approximation technique is the fact that it still needs to pass through the whole decision tree. Nonetheless, it can easily be parallelized, or randomly sampled HDMR may be used \cite{GenyuanGeorgopoulos2006_RS_HDMR}, 
or some reinforcement learning algorithm that aims at this problem can be applied. The fact that HDMR enables a fast approximate minimization would still be worthwhile. 

The author would like to express his gratitude to RNDr. Ond\v{r}ej Pangr\'{a}c, Ph.D., for inspiring discussion about  discrete optimization, to Irena Dvo\v{r}\'{a}kov\'{a}, prom. fil., for significant help with the language of the manuscript, and finally to Ing. V\'{a}clav \v{S}m\'{i}dl, Ph.D., for constructive criticism and encouragement.

\bibliographystyle{elsarticle-num}    


%
\end{document}